\documentclass[11pt]{amsart}

\usepackage[all]{xypic}

\usepackage{latexsym}

\usepackage{amssymb}
\usepackage{amsfonts}
\usepackage{amscd}
\usepackage{amsmath,amsthm}

\newtheorem{lemma}{Lemma}[section]

\newtheorem{lem}[lemma]{Lemma}

\newtheorem{thm}[lemma]{Theorem}
\newtheorem{cor}[lemma]{Corollary}

{

}

\theoremstyle{definition}

\theoremstyle{remark}

\numberwithin{equation}{section}

\newenvironment{pf}{\noindent{\bf Proof.}}{\hfill $\square$\medskip}

\def\CC{{\mathbb C}}

\def\NN{{\mathbb N}}

\def\PP{{\mathbb P}}

\def\ZZ{{\mathbb Z}}

\def\aol{{\bar a}}

\def\eol{{\bar e}}

\def\xol{{\bar x}}
\def\yol{{\bar y}}

\def\Qol{{\bar Q}}

\def\0ol{{\bar 0}}
\def\1ol{{\bar 1}}
\def\2ol{{\bar 2}}
\def\ol2{{\bar 2}}
\def\3ol{{\bar 3}}
\def\4ol{{\bar 4}}
\def\5ol{{\bar 5}}
\def\6ol{{\bar 6}}
\def\7ol{{\bar 7}}
\def\8ol{{\bar 8}}
\def\9ol{{\bar 9}}

\def\bold0{{\bf 0}}
\def\bold1{{\bf 1}}
\def\bold2{{\bf 2}} 
\def\bold3{{\bf  3}}
\def\bold4{{\bf 4}}
\def\bold5{{\bf 5}}
\def\bold6{{\bf 6}}
\def\bold7{{\bf 7}}
\def\bold8{{\bf 8}}
\def\bold9{{\bf 9}}

\def\P2Skly{\PP^2_{Skly}}

\def\coker{\operatorname {coker}}

\def\End{\operatorname {End}}

\def\Hilb{\operatorname {Hilb}}

\def\ker{\operatorname {ker}}

\def\th{\operatorname {th}}    
\def\Tor{\operatorname {Tor}}

\def\dim{\operatorname{dim}}

\def\End{\operatorname{End}}

\def\Fdim{{\sf Fdim}}

\def\Gr{{\sf Gr}}

\def\id{\operatorname{id}}

\def\liminj{\varinjlim}

\def\Mod{{\sf Mod}}

\def\Proj{\operatorname{Proj}}

\def\QGr{\operatorname{\sf QGr}}

\def\ul1{\operatorname{\underline{1}}}

\def\l{\leftarrow}

\def\l{\lambda}

\def\s{\sigma}

\def\fA{{\mathfrak A}}

\def\sB{{\sf B}}
\def\sC{{\sf C}}

\def\sK{{\sf K}}

\def\sX{{\sf X}}

\def\cal{\mathcal}

\def\cH{{\cal H}}

\def\cO{{\cal O}}

\def\Qcoh{{\sf Qcoh}}

\def\dirlim{\mathop{\vtop{\baselineskip -100pt\lineskip -1pt\lineskiplimit 0pt
\setbox0\hbox{lim}\copy0\hbox to \wd0{\rightarrowfill}}}\limits}
\def\invlim{\mathop{\vtop{\baselineskip -100pt\lineskip -1pt\lineskiplimit 0pt
\setbox0\hbox{lim}\copy0\hbox to \wd0{\leftarrowfill}}}\limits}

\def\I11{{1 \kern -0.8pt \! \mbox{l}}}
\def\mumu{{\mu\kern-4.2pt\mu}}
\def\bfmu{{\mu\kern-4.2pt\mu}}
\def\2slash{\backslash \! \backslash}

\def\boxtimes{\setbox0\hbox{$\Box$}\copy0\kern-\wd0\hbox{$\times$}}

\date{}                                           

\begin{document}

\title[Shift equivalence and modules over path algebras]{Shift equivalence and a category equivalence involving graded modules over path algebras of quivers}

\author{S. Paul Smith}

\address{ Department of Mathematics, Box 354350, Univ.
Washington, Seattle, WA 98195}

\email{smith@math.washington.edu}

\thanks{S. P. Smith was partially supported by NSF grant DMS 0602347}


\begin{abstract}
In this paper we associate an abelian category to a finite directed graph and prove the categories arising from two graphs are equivalent if the incidence matrices of the graphs are shift equivalent. The abelian category is the quotient  of the category of graded vector space representations of the quiver obtained by making the graded representations that are the sum of their finite dimensional submodules isomorphic to zero.

Actually, the main result in this paper is that the abelian categories are equivalent if the incidence matrices are strong shift equivalent. That result is combined with an earlier result of the author to prove that if the incidence matrices are shift equivalent, then the associated abelian categories are equivalent.

Given William's Theorem that subshifts of finite type associated to two directed graphs are conjugate if and only if the graphs are strong shift equivalent, our main result can be reformulated as follows: if the subshifts associated to two directed graphs are conjugate, then the categories associated to those graphs are equivalent.

\end{abstract}

\maketitle

\pagenumbering{arabic}


\setcounter{section}{0}

\section{Introduction}

\subsection{}
This paper proves a result of the following general type:  if two graphs are equivalent in 
an appropriate  sense, then certain algebraic objects associated to them are equivalent in a 
corresponding sense.  
We associate an abelian category to a directed graph and prove the categories arising from two graphs are equivalent if the graphs are equivalent in an appropriate  sense.

A nice paper by Bates and Pask \cite{BP}  contains results of this general type.
They prove various isomorphisms and Morita equivalences between graph 
C$^*$-algebras that unify earlier results for graph C$^*$-algebras and Cuntz-Kreiger algebras (references for some of those earlier results can be found in the opening paragraph of \cite{BP}). Raeburn's monograph \cite{R} 
is a comprehensive  treatment of graph C$^*$-algebras.
 Analogous, but algebraic as opposed to C$^*$-algebraic, results for Leavitt path algebras have been proved by Abrams et al. See, for example, \cite{ALPS} and \cite{AT} and the references therein. 
 
\subsection{}
The graph equivalences alluded to in the opening paragraph of this introduction 
arise in the theory of subshifts of finite type. 

A shift space $(X,\s)$ over a finite alphabet $\fA$ is a compact subset $X$ of $\fA^\ZZ$ that is stable under the shift map  $\sigma$ defined by $\s(f)(n)=f(n+1)$. 

To a directed graph $Q$ with finitely many edges and no sources or sinks one may associate a shift space $X_Q$, called the {\it edge shift} of $Q$, whose alphabet is the set of arrows in $Q$
(see \cite[Defn. 2.2.5]{LM}). Edge shifts are  subshifts of finite type. 

One may also associate a subshift of finite type $X_A$ to every square 0-1 matrix $A$
having no zero rows or columns  (see \cite[Defn. 2.3.7]{LM}).  The alphabet 
is now the set of vertices for the directed graph whose incidence matrix is $A$. 

The equivalence between shift spaces that concerns us is {\it conjugacy}.
Shift	spaces $(X,\s_X)$ and $(Y, \s_Y)$	are {\sf conjugate}, or {\sf topologically conjugate}, 
denoted $X \cong Y$, if there is a
homeomorphism $\phi:X \to Y$ such that $\s_Y \circ \phi = \phi \circ \s_X$.

By \cite[Prop. 2.3.9]{LM}, every subshift of finite type is conjugate to an edge shift $X_Q$ for some directed graph $Q$. 

One may associate to $Q$ a new graph $Q^{[2]}$ (see \cite[Defn. 2.3.4]{LM}) 
whose vertices are the edges of $Q$
and $Q^{[2]}$ has an arrow from $e$ to $e'$ if $e$ terminates where $e'$ begins. The incidence matrix for $Q^{[2]}$, which we denote by $B_Q$,   
is a 0-1 matrix such that $X_Q \cong X_{B_Q}$. Thus every subshift of finite type is conjugate to a shift described by a 0-1 matrix. 
Conversely, every shift described by a 0-1 matrix $A$ is conjugate to an edge shift: if $Q^A$ is the directed graph with incidence matrix $A$, then $X_{Q^A}$ is conjugate to $X_Q$  \cite[Exer. 1.5.6 and Prop. 2.3.9]{LM}. Hence subshifts of finite type are edge shifts or shifts associated to 0-1 matrices.

\subsection{The results}
\label{sect.1.1}
Throughout $k$ is a field and $Q$ a directed graph, or quiver, with a finite number of vertices and arrows---loops and multiple arrows between vertices are allowed. 

We write $kQ$ for the path algebra of $Q$. The finite paths, including the trivial path at each vertex, in $Q$
form a basis for $kQ$ and multiplication is given by concatenation of paths. 

We adopt the convention that the incidence matrix of $Q$ is $C=(c_{ij})$ where 
$$
c_{ij} = \hbox{the number of arrows from $j$ to $i$.}
$$
Given a square $\NN$-valued matrix we write $Q^C$ for the directed graph with incidence matrix $C$.

We make $kQ$ an $\NN$-graded algebra by declaring that a path is homogeneous of degree equal to its length.  The category of $\ZZ$-graded left $kQ$-modules with degree-preserving homomorphisms is denoted by $\Gr kQ$ and we write $\Fdim kQ$ for its full subcategory of consisting of modules that are the sum of their
finite-dimensional submodules. Since $\Fdim kQ$ is a Serre subcategory of $\Gr kQ$ (it is, in fact, a localizing subcategory) we may form the 
quotient category 
$$
\QGr kQ:=\frac{\Gr kQ}{\Fdim  kQ}.
$$

The main results in this paper is the following theorem and its consequences. 

\begin{thm}
\label{thm1}
Let $L$ and $R$ be $\NN$-valued matrices such that $LR$ and $RL$ make sense.
Let $Q^{LR}$ be the quiver with incidence matrix $LR$ and $Q^{RL}$ the quiver with incidence matrix $RL$.
There is an equivalence of categories
$$
\QGr kQ^{LR} \equiv \QGr kQ^{RL}.
$$
\end{thm}

Strong shift equivalence (see section \ref{sect.sse} for its definition) is an equivalence relation on square matrices with entries in $\NN$ that is important  in symbolic dynamics (see section \ref{sect.sse} below).  By interpreting a square matrix with entries in $\NN$ 
as an incidence matrix, an equivalence relation on square matrices with entries in $\NN$ is the same thing as an equivalence relation on finite directed graphs.

\begin{thm}
\label{thm4}
If the incidence matrices for $Q$ and $Q'$ are strong shift equivalent, then $\QGr(kQ) \equiv \QGr(kQ')$.
\end{thm}

Given William's Theorem (see Theorem \ref{thm.W} below), Theorem \ref{thm4} can be restated as follows:
 Let $(X,\s)$ and $(X,\s')$ be subshifts of finite type, and $Q$ and $Q'$ directed graphs such that 
$(X,\s) = X_Q$ and $(X,\s')=X_{Q'}$. If  $(X,\s)$ and $(X,\s')$ are conjugate, then the categories 
 $\QGr(kQ)$ and $\QGr(kQ')$ are equivalent.
 
It is difficult to decide if two given matrices are strong shift equivalent. It is not known whether the strong shift equivalence problem is decidable. However, there is a weaker notion, shift equivalence (see section 
\ref{sect.se} for its definition), and Kim and Roush \cite{KR1} have shown that the shift equivalence problem is decidable. Strong shift equivalence implies shift equivalence but the question of whether the two notions were
the same was open for over twenty years before Kim and Roush \cite{KR2}  gave an example in 1999 showing  shift equivalence does not imply strong shift equivalence. 

If two incidence matrices $A$ and $B$ are shift equivalent, then $A^\ell$ is strong shift equivalent to $B^\ell$
for some integer $\ell$. Theorem 1.8 in \cite{Sm2} says that if $Q^{(\ell)}$ is the directed graph whose incidence matrix is the $\ell^{\th}$ power of the incidence matrix for $Q$, then $\QGr kQ$ is equivalent to $\QGr (kQ^{(\ell)})$.\footnote{In symbolic dynamics $Q^{(\ell)}$ is called the $\ell^{\th}$ higher power graph of $Q$ \cite[Defn. 2.3.10]{LM}.}
Combining \cite[Thm. 1.8]{Sm2} with Theorem \ref{thm4} gives the following.

\begin{cor}
\label{cor0}
If the incidence matrices for $Q$ and $Q'$ are shift equivalent, then $\QGr(kQ) \equiv \QGr(kQ')$.
\end{cor}

Given $Q$, define $Q^{[n]}$ to be the following quiver: its vertices are the paths of length $n$ in $Q$;
if $p$ and $q$ are paths of length $n$ in $Q$ 
there is an arrow in $Q^{[n]}$ from $p$ to $q$ if there is a path of length $n+1$ in $Q$ that begins with $p$ and ends with $q$.

\begin{cor}
\label{cor2}
For all integers $n \ge 2$, $\QGr (kQ) \equiv \QGr (kQ^{[n]})$. 
\end{cor}

\medskip
{\bf Acknowledgements.} 
I wish to thank Doug Lind for introducing me to the notion of shift equivalence and for 
useful discussions about symbolic dynamics and related matters.


\section{(Strong) shift equivalence}

\subsection{Strong shift equivalence and Williams's Theorem}
\label{sect.sse}

Let $A$ and $B$ be square matrices with  entries in $\NN$.
An {\sf elementary strong shift equivalence} between $A$ and $B$ is a pair of matrices $L$ and $R$
with non-negative integer entries such that
$$
A=LR  \qquad \hbox{and} \qquad B=RL.
$$
We say $A$ and $B$ are {\sf strong shift equivalent} 
 if there is a chain of elementary strong shift equivalences from $A$ to $B$.

The following fundamental result explains the importance of strong shift equivalence.  
 
\begin{thm}
[Williams]
\cite[Thm. A]{W}
\label{thm.W}
Let $A$ and $B$ be square $\NN$-valued matrices and $X_A$ and $X_B$ the associated subshifts of finite 
type. Then $X_A \cong X_B$ if and only if $A$ and $B$ are strong shift equivalent.
\end{thm}

A proof of Theorem \ref{thm.W} can also be found at \cite[Thm. 7.2.7]{LM}.

\subsection{}
 
The matrices 
$$
A=\begin{pmatrix} 1 & 1 \\ 1 & 1 \end{pmatrix}=
\begin{pmatrix} 1 \\ 1 \end{pmatrix}\begin{pmatrix} 1 & 1 \end{pmatrix}
\qquad \hbox{and} \qquad 
B=(2)=\begin{pmatrix} 1 & 1 \end{pmatrix}\begin{pmatrix} 1 \\ 1 \end{pmatrix}
$$
are  strong shift equivalent. The corresponding quivers are
\vskip .1in
 $$
 \begin{array}{cc}
Q^A= \qquad  \UseComputerModernTips
\xymatrix{ 
\bullet     \ar@/^/[rr]  \ar@(ul,dl)[]  && 
 \bullet \ar@/^/[ll]  \ar@(ur,dr)[] 
}
&
\phantom{xxxxxxxxx}
Q^B= \qquad \UseComputerModernTips
\xymatrix{
\bullet    \ar@(ul,dl)[]    \ar@(ur,dr)[] 
}
\end{array}
$$

\vskip .1in

\noindent
By Theorem \ref{thm1}, $\QGr kQ^A \equiv \QGr kQ^B$.

By \cite{Sm2}, there are ultramatricial $k$-algebras $S(Q^A)$ and $S(Q^B)$ such that $\QGr kQ^A \equiv
\Mod S(Q^A)$ and $\QGr kQ^B \equiv \Mod S(Q^B)$. The Bratteli diagram for $S(Q^A)$ is 
$$
\UseComputerModernTips
\xymatrix{ 
1 \ar[r] \ar[dr] & 2 \ar[r] \ar[dr] & 4 \ar[r] \ar[dr] &  8 \ar[r] \ar[dr] &  16 \ar[r] \ar[dr] &  \cdots
\\
1 \ar[r] \ar[ur] & 2 \ar[r] \ar[ur] & 4 \ar[r] \ar[ur] &  8 \ar[r] \ar[ur] &   16   \ar[r] \ar[ur] &  \cdots
}
$$
and that for $S(Q^B)$ is 
$$
\UseComputerModernTips
\xymatrix{ 
 1 \ar@<.5ex>[r] \ar@<-.5ex>[r]   & 2 \ar@<.5ex>[r] \ar@<-.5ex>[r]  & 4 \ar@<.5ex>[r] \ar@<-.5ex>[r]  &  8 \ar@<.5ex>[r] \ar@<-.5ex>[r]  &
16  \ar@<.5ex>[r] \ar@<-.5ex>[r]  &  \cdots
}
$$

\subsection{}
No general procedure is known to decide if two matrices are strong shift equivalent. 
The shortest known sequence of elementary strong shift equivalences proving that  
$$
\begin{pmatrix}
  1    &  3  \\
  2    &  1
\end{pmatrix}
\qquad \hbox{and} \qquad 
\begin{pmatrix}
  1    &  6  \\
  1    &  1
\end{pmatrix}
$$
are strong shift equivalent was found by a computer search  \cite[Ex. 7.3.12]{LM} and en route from the first to the second matrix one passes through the incidence matrix  for the graph
\vskip .2in
   $$
Q=  \qquad \UseComputerModernTips
\xymatrix{
  && \ar[dll] \ar@/_2pc/[dll]  \bullet \ar@(ul,ur)[]   \ar[rrrr] \ar[drr] \ar@/^1pc/[drr]  &&&& 
\bullet \ar@/_2pc/[llll] \ar@/^1pc/[dll] 
\\
\bullet \ar@(dl,ul)[] \ar@/_1pc/[urr]  &&&&  \bullet  \ar@/^1pc/[ull]  \ar[urr]  \ar@/^2pc/[llll] \ar@/^1pc/[llll] &
}
$$
\vskip .4in

\noindent
Thus Theorem \ref{thm4} shows that 
$$
\QGr(kQ) \equiv \QGr(kQ') \equiv \QGr(kQ'')
$$
where
\vskip .2in
   $$
Q'= \qquad  \UseComputerModernTips
\xymatrix{
\bullet \ar@(dl,ul)[] \ar[rr] \ar@/_1pc/[rr]   \ar@/_2pc/[rr]   && \ar@/_2pc/[ll] \ar@/_1pc/[ll]  \bullet \ar@(ur,dr)  
}
\qquad \hbox{and} \qquad 
Q''= \qquad \UseComputerModernTips
\xymatrix{
\bullet \ar@(dl,ul)[]  \ar@/_1pc/[rr]  \ar[rr] \ar@/^1pc/[rr]   \ar@/^2pc/[rr]    \ar@/^3pc/[rr]    \ar@/^4pc/[rr] 
&& \ar@/^2pc/[ll]  \bullet \ar@(ur,dr)  
}
$$
\vskip .3in

\subsection{Shift equivalence}
\label{sect.se}

Two square matrices  $A$ and $B$ with non-negative integer entries are {\sf shift equivalent} if there is a positive integer $\ell$ and matrices $L$ and $R$ with non-negative integer entries such that
$$
AL=LB, \quad RA=RB, \quad A^{\ell}=LR,  \quad \hbox{and} \quad B^{\ell}=RL.
$$


\section{Proof of Theorem \ref{thm1}}
\label{sect.LR}

\subsection{Notation for quivers and path algebras}

Let $k$ be a field and $Q$ a finite quiver, i.e., a finite directed graph. We write $Q_0$ for its set of vertices 
and $Q_1$ for its set of arrows. If the arrow $a$ ends where the arrow $b$ starts we write $ba$ for the 
path ``first traverse $a$ then traverse $b$''. We write $Q_n$ for the set of paths of length $n$ in $Q$.

If $p$ is a path we write $s(p)$ for the vertex at which it starts and $t(p)$ for the vertex at which it terminates.

The path algebra $kQ$ has a basis given by the set of all finite paths, including the empty path and the trivial paths at each vertex. The multiplication in $kQ$ is the linear extension of that given by concatenation of paths, i.e.,
$$
q \times p =
\begin{cases} 
	qp & \text{if $s(q)=t(p)$}
	\\
	0 & \text{if $s(q)\ne t(p)$.}
\end{cases}
$$	

The algebra $kQ$ is $\NN$-graded with degree $n$ component equal to $kQ_n$, the linear span of the paths of length $n$. The subalgebra $kQ_0$ of $kQ$ is isomorphic to a product of $|Q_0|$
copies of $k$ and is therefore a semisimple ring. Each $kQ_n$ is a $kQ_0$-bimodule. The
multiplication in $kQ$ gives an isomorphism
$$
 (kQ_1)^{\otimes n} \cong kQ_n
$$
of $kQ_0$-bimodules where the tensor product on the left-hand side is taken over $kQ_0$. 
It follows that $kQ$ is isomorphic to the tensor algebra over $kQ_0$ of $kQ_1$,
$$
kQ \cong T_{kI}(kQ_1).
$$

\subsection{}
Let $i$ and $j$ be positive integers.

Let $kI$ denote the ring of $k$-valued functions on $[i]=\{1,\ldots,i\}$ with pointwise addition and multiplication. Similarly, $kJ$ denotes the ring of $k$-valued functions on $[j]=\{1,\ldots,j\}$. 
We identify $kI \otimes_k kJ$ with  the ring of $k$-valued functions on the Cartesian product $[i] \times [j]$. 
The category of $kI$-$kJ$-bimodules  is equivalent to the category of $kI \otimes_k kJ$-modules 
and we write $E_{pq}$ for the simple $kI$-$kJ$-bimodule corresponding to $(p,q) \in [i] \times [j]$; i.e., $E_{pq}$ is a copy of  $k$ supported at $(p,q)$. 

\subsection{}
Let $L=(\ell_{pq})$ be an $i \times j$ matrix over $\NN$ and $R=(r_{st})$ a $j \times i$ matrix over $\NN$. Let
$Q^{LR}$ be the directed graph with incidence matrix $LR$ and $Q^{RL}$ the directed graph 
with incidence matrix $RL$.

Since it is unlikely to cause confusion we will also use the letter $L$ to denote the $kI$-$kJ$-bimodule
$$
L:=\bigoplus_{\substack{p \in [i] \\ q \in [j]}} (E_{pq})^{\oplus \ell_{pq}}.
$$
In a similar way we define the $kJ$-$kI$-bimodule
$$
R:=\bigoplus_{\substack{t \in [i] \\ s \in [j]}} (E_{st})^{\oplus r_{st}}.
$$

The linear span in $kQ^{LR}$ of the arrows in $Q^{LR}$ is isomorphic
to $L \otimes_{kJ} R$ as a $kI$-$kI$-bimodule. 
We identify the path algebras $kQ$ and $kQ'$ with the following tensor algebras:
$$
kQ^{LR}=T_{kI}(L \otimes_{kJ} R) = \bigoplus_{n=0}^\infty (L \otimes_{kJ} R)^{\otimes n}
$$
and
$$
kQ^{RL}=T_{kJ}(R \otimes_{kI} L) = \bigoplus_{n=0}^\infty (R \otimes_{kI}L)^{\otimes n}.
$$
We give $kQ^{LR}$ its standard grading by declaring that $kI$ is its degree-zero component and 
$L \otimes_{kJ} R$ its degree-one component.

\subsection{}
Since $kQ^{LR}$ is the tensor algebra of the $kI$-bimodule $L \otimes_{kJ} R$, a graded
 left $kQ^{LR}$-module is a pair $(M,\l)$ consisting of a graded left $kI$-module $M$ and a  homomorphism 
$$
\l: L \otimes_{kJ} R \otimes_{kI} M \to M
$$
of left $kI$-modules such that $\l(L \otimes_{kJ} R \otimes_{kI} M_n) \subset M_{n+1}$ for all $n$. A 
homomorphism $(M,\l) \to (M',\l')$ of graded $kQ^{LR}$-modules is a homomorphism $\theta:M \to
M'$ of graded $kI$-modules such that
$$
\theta \circ \l = \l' \circ (\id_L \otimes \id_R \otimes \theta).
$$

\subsection{}
We now define functors 
$$
\UseComputerModernTips
\xymatrix{
\Gr \big( kQ^{LR} \big)  \ar@/^1pc/[rrr]^F & && \Gr \big( kQ^{RL} \big)\ar@/^1pc/[lll]^{F'}.
}
$$
If $M$ is a graded left $kQ^{LR}$-module we define
$$
F(M,\l):= (R \otimes_{kI} M, \id_R \otimes \l)
$$
with the grading
$$
(R \otimes_{kI}M)_n:= R \otimes_{kI}M_{n}.
$$
The action of the degree-one component of $kQ^{RL}$, which is $R \otimes_{kI} L$, on 
the dgree $n$ component $F(M,\l)_n$ is
\begin{align*}
(\id_R \otimes \l)\big((R \otimes_{kI} L) \otimes_{kJ} (R \otimes_{kI} M)_n\big) 
& = R \otimes_{kI} \l(L \otimes_{kJ} R \otimes_{kI} M_{n})
\\
& \subset R \otimes_{kI} M_{n+1}
\\
&= ( R \otimes_{kI} M)_{n+1}
\end{align*}
so $R \otimes_{kI} M$ really is a {\it graded} $kQ'$-module.  

If $\theta:(M,\l) \to (M,\l')$ is a homomorphism of graded $kQ^{LR}$-modules we define 
$$
F(\theta):=\id_R \otimes \theta.
$$
It is easy to check that $F\theta:F(M,\l) \to F(M',\l')$ is a homomorphism of graded 
$kQ^{RL}$-modules. Hence $F$ is a functor.

The functor $F'$ is defined in a similar way.

Since $kI$ and $kJ$ are semisimple rings $F$ and $F'$ are exact functors.

\begin{thm}
\label{thm.LR}
Let $L$ be an $i \times j$ matrix over $\NN$ and $R$ a $j \times i$ matrix over $\NN$.  
Then the functors $F$ and $F'$ induce mutually quasi-inverse equivalences of categories 
$$
\QGr\big(kQ^{LR}\big) \equiv \QGr\big(kQ^{RL}\big).
$$
\end{thm}
\begin{pf}
Let $(M,\l) \in \Gr \big(kQ^{LR}\big)$. Then
\begin{align*}
F'F(M,\l) & =F'(R \otimes_{kI} M,\id_R \otimes \l)
\\
& = (L \otimes_{kJ} R \otimes_{kI} M,\id_L \otimes \id_R \otimes \l).
\end{align*}
We define $\tau_M:F'FM \to M$ by $\tau_M(x \otimes y \otimes m):=\l(x \otimes y \otimes m)$, 
Since $\tau_M=\l$ it is a tautology that
$$
\tau_M \circ (\id_L \otimes \id_R \otimes \l) = \l \circ (\id_L \otimes \id_R \otimes \tau_M)
$$
whence $\tau_M$ is a homomorphism of $kQ^{LR}$-modules.
Since $$(F'FM)_n = L \otimes_{kJ} (FM)_{n-1}=   L \otimes_{kJ} R \otimes_{kI} M_{n-1}$$ we have
$$
\tau_M\big((F'FM)_n\big) =\l(L \otimes_{kJ} R \otimes_{kI} M_{n-1}) \subset M_{n}.
$$
Hence $\tau_M$ is a homomorphism of graded $kQ^{LR}$-modules. 

The above shows that $$\tau:F'F \to \id_{\Gr ( kQ^{LR})}$$ is a natural transformation.

Since $F$ and $F'$ are exact functors that send finite dimensional modules to finite dimensional modules they induce functors between the quotient categories $\QGr ( kQ^{LR})$ and $\QGr ( kQ^{RL})$, say $f$ and $f'$. It follows that $\tau$ induces a natural transformation from $f'f$ to $\id_{\QGr ( kQ^{LR})}$. 
We will now show this induced natural transformation is an isomorphism of functors. A similar argument will show $ff'$ is isomorphic to $\id_{\QGr ( kQ^{RL})}$. The proof of the theorem will then be complete.

Write $V=L \otimes_{kJ} R$. 

\underline{Claim}: If $M \in \Gr ( kQ^{LR})$, then $F'FM \cong (kQ^{LR})_{\ge 1} \otimes_{kQ^{LR}} M$ as graded left $kQ^{LR}$-modules.
\underline{Proof}:
Let $(M,\l)$ be a graded left $kQ^{LR}$-module. Then 
$$
F'F(M,\l)  = (V \otimes_{kI} M, \id_V \otimes \l).
$$
We make the identification $(kQ^{LR})_{\ge 1}=V \otimes_{kI} kQ^{LR}$ so the formula
$$
 \theta(v \otimes a \otimes m) := v \otimes am
 $$
 for $v \in V$, $a \in kQ^{LR}$, and $m \in M$, defines an isomorphism of left $kI$-modules
$$
\theta:  (kQ^{LR})_{\ge 1} \otimes_{kQ^{LR}} M = V \otimes_{kI} kQ^{LR}  \otimes_{kQ^{LR}} M \, \longrightarrow \,
  V \otimes_{kI}   M.
$$
If $v' \in V$ and 
$v \otimes a \otimes m \in V \otimes_{kI} kQ^{LR}  \otimes_{kQ^{LR}} M$, then
\begin{align*}
\theta\big(v'.(v \otimes a \otimes m)\big) &= \theta(v'  \otimes va \otimes m) 
\\
& = v'  \otimes va m
\\
& = (\id_V \otimes \l)(v' \otimes v \otimes am)
\\
& = v'.(v  \otimes am) 
\\
& = v' .\theta(v \otimes a \otimes m) 
\end{align*}
so $\theta$ is a homomorphism, and therefore an isomorphism, of left $kQ^{LR}$-modules. In fact, 
$$
\theta:(kQ^{LR})_{\ge 1}\otimes_{kQ^{LR}} M \to F'F(M,\l)
$$
is an isomorphism of {\it graded} $kQ$-modules because if $v \otimes a \otimes m$ is a homogeneous
element of $V \otimes_{kI} kQ \otimes_{kQ} M =  (kQ)_{\ge 1}\otimes_{kQ} M$, then
$$
\deg(v \otimes a \otimes m)=1+\deg (am);
$$
however, $F'F(M,\l)_n = V \otimes_{kI} M_{n-1}$ so, as an element of $F'F(M,\l)$, $\deg (v \otimes am)
= 1+\deg(am)$ so
$$
\deg \theta(v \otimes a \otimes m) = \deg (v \otimes am) =  \deg (v \otimes a \otimes m);
$$
i.e., $\theta$ is a degree-preserving map so an isomorphism in $\Gr kQ$. 
This completes the proof of the claim. $\lozenge$

The claim shows that the homomorphisms 
\begin{align*}
\eta_M:F'F(M,\l)  & = (V \otimes_{kI} M, \id_V \otimes \l) \to (kQ)_{\ge 1} \otimes_{kQ} M
\\
\eta_M(v \otimes m) &=v \otimes m
\end{align*}
produce an isomorphism of functors
$$
\eta:F'F  \to (kQ)_{\ge 1} \otimes_{kQ}-.
$$

Consider the diagram
 $$
  \UseComputerModernTips
\xymatrix{
&&  F'FM \ar[r]^{\tau_M} \ar[d]_{\eta_M} & M \ar@{=}[d]
\\
0 \ar[r] & \Tor^{kQ}_1(kI,M) \ar[r] &  (kQ)_{\ge 1} \otimes_{kQ} M \ar[r]_<<<<<{\mu}& M \ar[r]& kI \otimes_{kQ} M \ar[r] & 0
}
$$
where the bottom row is the exact sequence obtained by applying $-\otimes_{kQ} M$ to the exact sequence of $kQ$-bimodules $0 \to (kQ)_{\ge 1} \to kQ \to kI \to 0$ and $\mu$ is the multiplication in $kQ$.
Since $\tau_M(v \otimes m)=\l(v \otimes m)=vm$, the square commutes. Hence
$$
\ker \tau_M \cong \Tor^{kQ}_1(kI,M) 
\qquad \hbox{and} \qquad 
\coker \tau_M \cong kI \otimes_{kQ} M.
$$
Both these modules are annihilated by $(kQ)_{\ge 1}$ so belong to $\Fdim(kQ)$. Therefore, after passing 
to $\QGr kQ$, the diagram yields a commutative square in which $\mu$ and $\eta_M$ become isomorphisms.
It follows that $\tau_M$ becomes an isomorphism in $\QGr kQ$, and hence that $\tau:F'F \to \id_{\QGr kQ}$ is an isomorphism of functors as claimed.

Given the symmetry of the situation we can reverse the roles of $L$ and $R$ and repeat the previous argument to produce an isomorphism of functors $FF' \to \id_{\QGr kQ'}$. This completes the proof that $\QGr kQ$ is equivalent to $\QGr kQ$.
\end{pf}


\section{Proof of Theorem \ref{thm4}}

Suppose the $\NN$-valued matrices $A$ and $B$ are strong shift equivalent.
By definition, there is a sequence of matrices
$$
A=A_1, \, A_2, \, \ldots, \, A_n=B
$$
 and elementary strong shift equivalences between $A_i$ and $A_{i+1}$ for $1 \le i \le n-1$.
If $Q^{A_i}$ is the quiver with incidence matrix $A_i$, then repeated applications of Theorem \ref{thm1} 
show that
$$
\QGr kQ^{A_1} \equiv  \QGr kQ^{A_2} \equiv  \cdots \equiv \QGr kQ^{A_n}
$$
thereby proving Theorem \ref{thm4}.


\subsection{In-splitting and out-splitting}
I am grateful to Min Wu for telling me that Theorem \ref{thm1} applies to in-splittings and out-splittings.

Lind and Marcus define and discuss in-splittings and out-splittings of a directed graph in section 2.4 of \cite{LM}. As the name suggests, in-splitting involves replacing one vertex $v$ by several, say $n$, 
vertices $v_1,\ldots,v_n$ and replacing each arrow $a$ ending at $v$ by $n$ arrows $a_1,\ldots,a_n$ 
where $a_i$ starts where $a$ does and ends at $v_i$. Actually, in-splitting is a more general process than this, but the basic idea is along the lines just described. Out-splitting is an analogous process, now based on the arrows leaving a vertex.

The important point for us is that if $Q'$ is obtained from $Q$ by an in-splitting or an out-splitting there is an elementary strong shift equivalence between their incidence matrices (see \cite[Thm. 2.4.12]{LM} and \cite[Exer. 2.4.9]{LM}). Theorem \ref{thm1} therefore yields the following result.

\begin{cor}
\label{cor1}
If $Q'$ is obtained from $Q$ by an in-splitting or out-splitting, then $\QGr kQ \equiv \QGr kQ'$.
\end{cor}


\section{Proof of Corollary \ref{cor2}}

\subsection{}
Let $Q$ be a finite quiver. Define $Q^{[2]}$ by setting
$$
Q^{[2]}_0:= Q_1
\qquad \hbox{and} \qquad Q^{[2]}_1:=Q_2
$$
with a length-two path $ba$ in $Q$ being considered as an arrow in $Q^{[2]}$ from $a$ to $b$. 

If $a$ and $b$ are arrows in $Q$, there is at most one arrow in  $Q^{[2]}$ from $a$ to $b$ so the incidence 
matrix for $Q^{[2]}$ is a 0-1 matrix.

The following example  from \cite[Example 2.7]{R} and \cite[Example 1.4.2]{LM} illustrates the construction:
$$
Q= \quad  \UseComputerModernTips
\xymatrix{
1    \ar@/^/[rr]^v  \ar@(ul,dl)[]_w  && 2 \ar@/^/[ll]^u
}
\qquad \qquad
Q^{[2]}=  \UseComputerModernTips
\xymatrix{
& w  \ar[dr]^{wv}  \ar@(ul,ur)[]^{ww}
\\
u   \ar[ur]^{wu}  \ar@/^/[rr]^{vu}    && v \ar@/^/[ll]^{uv}
}
$$

Returning to the general case, define a $|Q_0| \times |Q_1|$ matrix $L$  by 
$$
L_{ia}= \begin{cases}
		1 & \text{if $t(a)=i$}
		\\
		0 & \text{if $t(a)\ne i$}
	\end{cases}
$$	
for $i \in Q_0$ and $a \in Q_1$, 
and a $|Q_1| \times |Q_0|$ matrix $R$  by 
$$
R_{ai}= \begin{cases}
		1 & \text{if $s(a)=i$}
		\\
		0 & \text{if $s(a)\ne i$.}
	\end{cases}
$$
Then $LR$ is the $|Q_0| \times |Q_0|$ matrix with entries
\begin{align*}
(LR)_{ij}  & = \big\vert \{ a \in Q_1 \; | \; t(a)=i \, \hbox{ and } \, s(a)= j\} \big\vert
\\
&= \hbox{the number of arrows in $Q$ from $j$ to $i$}
\end{align*}
and $RL$ is the $|Q_1| \times |Q_1|$ matrix with entries
\begin{align*}
(RL)_{ab}  & = \big\vert \{ i \in Q_0 \; | \; s(a)=i =t(b)\} \big\vert 
= \begin{cases}
		1 & \text{if $s(a)=t(b)$}
		\\
		0 & \text{if $s(a)\ne t(b)$.}
	\end{cases}
\end{align*}
Therefore $LR$ is the incidence matrix for $Q$ and $RL$ is the incidence matrix for $Q^{[2]}$. 
Theorem \ref{thm1} therefore gives an equivalence
\begin{equation}
\label{eq.step1}
\QGr kQ \equiv \QGr kQ^{[2]}.
\end{equation}

\subsection{}
Let $Q$ be a finite quiver and define $Q^{[n]}$ by setting
$$
Q^{[n]}_0:= Q_n
\qquad \hbox{and} \qquad Q^{[n]}_1:=Q_{n+1}
$$
with a path $a_n\ldots a_0$ of length $n+1$ in $Q$ being considered as an arrow in $Q^{[n]}$ from $a_{n-1}
\ldots a_0$ to $a_n\ldots  a_1$.

It is an easy exercise to show that
$$
\Big(Q^{[n-1]}\Big)^{[2]} = Q^{[n]}.
$$
Repeatedly applying the equivalence in (\ref{eq.step1}) gives a chain of equivalences
$$
\QGr kQ \equiv \QGr kQ^{[2]}  \equiv\cdots  \equiv \QGr kQ^{[n]}
$$
thereby proving Corollary \ref{cor2}.


\begin{thebibliography}{12}
  
  \bibitem{ALPS}
  G. Abrams, A. Louly, E. Pardo, and C. Smith, 
  Flow invariants in the classification of Leavitt
path algebras, {\it Journal of Algebra,} {\bf 333} (2011) 202-231.
  
  \bibitem{AT}
G.   Abrams and M. Tomforde,
 Isomorphism and Morita equivalence of graph algebras,
 {\it Trans. Amer. Math. Soc.,} {\bf 363} (2011) 3733-3767. 
 

 \bibitem{BP}
 T. Bates and D. Pask,
 Flow equivalence of graph algebras,
 {\it  Ergodic Theory and Dynamical Systems,} {\bf 24} (2004) 367-382. 

\bibitem{KR1}
 K. H. Kim and F. W. Roush, Decidability of shift equivalence, {\it Proceedings of Maryland Special Year in Dynamics 1986-87,}  Springer-Verlag Lecture Notes in Math. 1342 (1988), 374-424.
 
 \bibitem{KR2}
 K. H. Kim and F. W. Roush, Williams' conjecture is false for irreducible subshifts, Annals of Math, {\bf 149} (199)
 545-558.

\bibitem{LM}
D. Lind and B. Marcus, {\it An Introduction to Symbolic Dynamics and Coding,} 
Camb. Univ. Press, Cambridge, 1995.
 
 \bibitem{R}
 I. Raeburn, {\it Graph algebras,} CBMS Regional Conference Series in Mathematics, 103, Published for the Conference Board of the Mathematical Sciences, Washington, DC, by the Amer. Math. Soc., Providence, RI, 2005. MR 2135030 (2005k:46141)






\bibitem{Sm1}
S.P. Smith, 
The non-commutative scheme having a free algebra as a homogeneous coordinate ring,
arXiv:1104.3822

\bibitem{Sm2}
S.P. Smith,  
Category equivalences involving graded modules over path algebras of quivers,
arXiv:1107.3511 


\bibitem{W}
R.F. Williams, Classification of subshifts of finite type,
{\it Ann. of Math.,} {\bf 98} (1973) 120-153.


\end{thebibliography}
\end{document}